\documentclass{fic-l}
\usepackage{latexsym}
\usepackage{amssymb}
\usepackage{amsthm}
\usepackage{times}
\usepackage[all]{xy}
\usepackage{color}

\newtheorem{thm}{Theorem}[section]

\newtheorem{lem}[thm]{Lemma}
\newtheorem{prop}[thm]{Proposition}
\newtheorem{exa}[thm]{Example}
\newtheorem{defi}[thm]{Definition}

\newcommand{\Jac}{\mathrm{Jac}\,}
\newcommand{\End}{\mathrm{End}\,}

\newcommand{\Aut}{\mathrm{Aut}\,}
\newcommand{\GL}{\mathrm{GL}}

\newcommand{\ord}{\mathrm{ord}}
\renewcommand{\H}{({\mathbb Z}/\ell {\mathbb Z})}
\newcommand{\Gb}{({\mathbb Z}/\ell {\mathbb Z})^b \rtimes {\mathbb Z}/p {\mathbb Z}}
\newcommand{\Ga}{({\mathbb Z}/\ell {\mathbb Z})^a \rtimes {\mathbb Z}/p {\mathbb Z}}

\newcommand{\bba}{\mathbb{A}}

\newcommand{\bbp}{\mathbb{P}}
\newcommand{\bbq}{\mathbb{Q}}
\newcommand{\bbz}{\mathbb{Z}}

\begin{document}
\title{Semi-direct Galois covers of the affine line}

\author[Gruendken]{Linda Gruendken}
\address{Department of Mathematics,
University of Pennsylvania,
David Rittenhouse Lab, 
209 South 33rd Street, 
Philadelphia, PA 19104-6395}\email{lindagr@math.upenn.edu}

\author[Hall-Seelig]{Laura Hall-Seelig}
\address{Department of Mathematics and Statistics,
Lederle Graduate Research Tower,
University of Massachusetts,
Amherst, MA 01003-9305}\email{hall@math.umass.edu}

\author[Im]{Bo-Hae Im}
\address{
Department of Mathematics,
Chung-Ang University,
221, Heukseok-dong, Dongjak-gu,
Seoul 156-756\\ South Korea
}
\email{imbh@cau.ac.kr}

\author[Ozman]{Ekin Ozman}
\address{Department of Mathematics,
University of Wisconsin-Madison,
480 Lincoln Drive,
Madison, WI 53706}\email{ozman@math.wisc.edu}

\author[Pries]{Rachel Pries}
\address{
Department of Mathematics,
Colorado State University,
Fort Collins, CO 80523-1874, USA
}
\email{pries@math.colostate.edu}

\author[Stevenson]{Katherine Stevenson}
\address{       Department of Mathematics,
California State University, 
18111 Nordhoff St,
Northridge, CA 91330-8313
}
\email{katherine.stevenson@csun.edu}

\date{\today}

\thanks{This project was initiated at the workshop WIN Women in Numbers in November 2008.  The authors would like to thank the Banff International Research Station for hosting the workshop and the National Security Agency, the Fields Institute, the Pacific Institute for the Mathematical Sciences, Microsoft Research, and University of Calgary for their financial support.  Author Pries was partially supported by NSF grant 07-01303. Author Im was partially supported by the Korea Science and Engineering 
Foundation (KOSEF) grant (No.~R01-2007-000-10660-0) funded by the Korea 
government (MOST). The authors would also like to thank the referee for helpful comments}
\subjclass{14H30, 14E20, 14F40}

\begin{abstract}
Let $k$ be an algebraically closed field of characteristic $p>0$. 
Let $G$ be a semi-direct product of the form $\Gb$ where $b$ is a positive integer
and $\ell$ is a prime distinct from $p$. 
In this paper, we study Galois covers $\psi:Z \to \bbp^1_k$ ramified only over $\infty$ with Galois group $G$.
We find the minimal genus of a curve $Z$ which admits a covering map of this form
and we give an explicit formula for this genus in terms of $\ell$ and $p$. 
The minimal genus occurs when $b$ equals the order $a$ of $\ell$ modulo $b$ and
we also prove that the number of curves $Z$ of this minimal genus 
which admit such a covering map is at most $(p-1)/a$ when $p$ is odd.
\end{abstract}

\maketitle

\section{Introduction}
Let $k$ be an algebraically closed field of characteristic $p>0$.  
In sharp contrast with the situation in characteristic $0$, there exist Galois covers $\psi:Z \to \bbp^1_k$ 
ramified only over infinity.  
By Abhyankar's Conjecture \cite{ab}, proved by Raynaud and Harbater \cite{ray}, \cite{har}, 
a finite group $G$ occurs as the Galois group of such a cover $\psi$ if and only if $G$ is quasi-$p$, 
i.e., $G$ is generated by $p$-groups.  This result classifies all the finite quotients of the fundamental group 
$\pi_1(\bba^1_k)$.  It does not, however, determine the profinite group structure of $\pi_1(\bba^1_k)$ 
because this fundamental 
group is an infinitely generated profinite group.

There are many open questions about Galois covers $\psi:Z \to \bbp^1_k$ ramified only over infinity.
For example, given a finite quasi-$p$ group $G$, what is the smallest integer $g$ for which there 
exists a cover $\psi:Z \to \bbp^1_k$ ramified only over infinity with $Z$ of genus $g$?  
As another example, suppose $G$ and $H$ are two finite
quasi-$p$ groups such that $H$ is a quotient of $G$.  
Given an unramified Galois cover $\phi$ of $\bba^1_k$ with group $H$, 
under what situations can one dominate $\phi$ with an unramified Galois cover $\psi$ of $\bba^1_k$
with Galois group $G$?
Answering these questions will give progress towards understanding how the finite quotients of 
$\pi_1(\bba^1_k)$ fit together in an inverse system.
These questions are more tractible for quasi-$p$ groups that are $p$-groups since   
the maximal pro-$p$ quotient $\pi_1^p(\bba^1_k)$ is free (of infinite rank) \cite{Shaf}.  

In this paper, we study 
Galois covers $\psi:Z \to \bbp^1_k$ ramified only over $\infty$ 
whose Galois group is a semi-direct product of the form $\Gb$, where $\ell$ is a prime distinct from $p$.
Such a cover $\psi$ must be a composition $\psi=\phi \circ \omega$ where 
$\omega:Z \to Y$ is unramified and 
$\phi: Y \to \bbp^1_k$ is an Artin-Schreier cover ramified only over $\infty$.
The cover $\phi$ has an affine equation $y^p-y=f(x)$ for some $f(x) \in k[x]$ with degree $s$ prime-to-$p$.
The $\ell$-torsion $\Jac(Y)[\ell]$ of the Jacobian of $Y$ is isomorphic to $(\bbz/\ell \bbz)^{2g_Y}$.
When $f(x)=x^s$, we determine how an automorphism $\tau$ of $Y$ of order $p$ acts on $\Jac(Y)[\ell]$.
This allows us to construct a Galois cover $\psi_a:Z_a \to \bbp^1_k$ ramified only over $\infty$ which dominates 
$\phi$, such that the Galois group of $\psi_a$ is $\Ga$ where $a$ is the order of $\ell$ modulo $p$
(Section~\ref{exist}).
We prove that the genus of $Z_a$ is minimal among all natural numbers 
that occur as the genus of a curve $Z$ which admits a covering map 
$\psi:Z \to \bbp^1_k$ ramified only over $\infty$ with Galois group of the form $\Gb$.
We also prove that the number of curves $Z$ of this minimal genus 
which admit such a covering map is at most $(p-1)/a$ when $p$ is odd  
(Section~\ref{unique}).

\section{Quasi-$p$ semi-direct products} \label{property}

We recall which groups occur as Galois groups of covers of $\bbp^1_k$ ramified only over $\infty$.

\begin{defi}\label{quasi-p} A finite group is a quasi $p$-group if it is generated by all of its Sylow $p$-subgroups.
\end{defi}

It is well-known that there are other equivalent formulations of the quasi-$p$ property, such as the next result.

\begin{lem}\label{lem:quasi-p} 
A finite group is a quasi $p$-group if and only if it has no nontrivial quotient group
whose order is relatively prime to $p$.
\end{lem}

The importance of the quasi-$p$ property is that it characterizes which finite groups occur as Galois groups 
of unramified covers of the affine line.

\begin{thm} \label{Tabconj}
A finite group occurs as the Galois group of a Galois cover of the projective line $\bbp^1_k$ ramified 
only over infinity if and only if it is a quasi-$p$ group.
\end{thm}

\begin{proof}
This is a special case of Abhyankar's Conjecture \cite{ab} which was jointly proved by Harbater \cite{har} and Raynaud \cite{ray}.
\end{proof}

We now restrict our attention to groups $G$ that are semi-direct products of the form $\Gb$.
The semi-direct product action is determined by a homomorphism $\iota:\bbz/p \bbz \to \Aut((\bbz/\ell \bbz)^b)$.

\begin{lem} \label{Lgroup}
Suppose a quasi-$p$ group $G$ is a semi-direct product of the form $\Gb$ for a positive integer $b$.
\begin{enumerate}
\item Then $G$ is not a direct product.
\item Moreover, $b \geq \ord_p(\ell)$ where $\ord_p(\ell)$ is the order of $\ell$ modulo $p$.
\end{enumerate}
\end{lem}

\begin{proof} 
Part (1) is true since $(\bbz/\ell \bbz)^b$ cannot be a quotient of the quasi-$p$ group $G$.
For part (2), the structure of a semi-direct product $\Gb$ 
depends on a homomorphism $\iota: \bbz/ p \bbz \to \Aut(\bbz / \ell \bbz)^b$.
By part (1), $\iota$ is an inclusion.  
Thus $\Aut(\bbz / \ell \bbz)^b \simeq \GL_b\H$ has an element of order $p$.
Now $$|\GL_b\H| = (\ell^b-1)(\ell^b-\ell)\cdots (\ell^b-\ell^{b-1}).$$ 
Thus $\ell^\beta \equiv 1\bmod p$ for some positive integer $\beta \leq b$ which implies $b \geq \ord_p(\ell)$.
\end{proof}

\begin{lem} \label{LexistGb}
If $a = \ord_p(\ell)$, then there exists a semi-direct product of the form $\Ga$ which is quasi-$p$.
It is unique up to isomorphism.
\end{lem}

\begin{proof}
If $a = \ord_p(\ell)$, then there is an element of order $p$ in $\Aut((\bbz / \ell \bbz)^a)$
and so there is an injective homomorphism 
$\iota: \bbz /p \bbz \to \Aut((\bbz / \ell \bbz)^a)$.
Thus there exists a non-abelian semi-direct product $G$ of the form $\Ga$.
To show that $G$ is quasi-$p$, 
suppose $N$ is a normal subgroup of $G$ whose index is relatively prime to $p$.
Then $N$ contains an element $\tau$ of order $p$.
By \cite[5.4, Thm.\ 9]{df}, since $G$ is not a direct product and $(\bbz / \ell \bbz)^a$ is normal in $G$, 
the subgroup $\langle \tau \rangle$ is not normal in $G$.
Thus $\langle \tau \rangle$ is a proper subgroup of $N$.
It follows that $\ell$ divides $|N|$ and so $N$ contains an element $h$ of order $\ell$ by Cauchy's theorem.
Recall that $\Aut((\bbz / \ell \bbz)^\beta)$ contains no element of order $p$ for any positive integer
$\beta < a$.
Thus the group generated by the conjugates of $h$ under $\tau$ has order divisible by $\ell^a$.
Thus $N=G$ and $G$ has no non-trivial quotient group whose order is relatively prime to $p$.
By Lemma \ref{lem:quasi-p}, $G$ is quasi-$p$.

The uniqueness follows from \cite[Lemma 6.6]{PS:survey}.
\end{proof}

\section{Explicit construction of $\Ga$-Galois covers of $\bba^1_k$} \label{exist}

In this section, we give concrete examples of Galois covers $\psi:Z \to \bbp^1_k$ 
ramified only over $\infty$ with Galois group of the form $\Ga$.
To compute the genus of the covering curve $Z$, we will need to determine the higher ramification
groups of $\psi$.

\begin{defi}
Let $L/K$ be a Galois extension of function fields of curves 
with Galois group $G$ and let $P,P'$ be primes of $K$ and $L$ such that $P'|P$. 
Let $\nu_{P'}$ and ${\mathcal O}_{P'}$ be the corresponding valuation function and valuation ring for $P'$.
For any integer $i \geq -1$, the $i$th ramification group of $P'|P$ is 
$$I_i(P'|P)=\{\sigma \in G \ | \ \nu_{P'}(\sigma(z)-z) \geq i+1 ,\forall z \in {\mathcal O}_{P'}\}.$$ 
\end{defi}

\begin{lem}
\label{lem:genus}
Suppose $f(x) \in k[x]$ is a polynomial of degree $s$ for a positive integer $s$ prime to $p$.
Let $\phi: Y \to \mathbb{P}^1_k$ be the cover of curves corresponding to the field extension 
$$k(x) \hookrightarrow k(x)[y]/ (y^p - y - f(x)).$$   
\begin{enumerate}
\item Then $\phi: Y \to \mathbb{P}^1_k$ is a Galois cover with Galois group $\bbz/ p \bbz$
ramified only at the point $P_\infty$ over $\infty$.
\item The $i$th ramification group at $P_\infty$ satisfies 
\begin{displaymath}
I_i = \left\{ \begin{array}{ll}
\mathbb{Z}/p\mathbb{Z} & \textrm{if } i \leq s\\
0 & \textrm{if } i > s. \end{array}\right.
\end{displaymath}
\item The genus $g_Y$ of $Y$ is equal to
$$g_Y=(p-1)(s-1)/2.$$
\end{enumerate}
\end{lem}

\begin{proof}
For part (1), note that the extension $k(x) \hookrightarrow k(x)[y]/ (y^p - y - f(x))$ is
cyclic of degree $p$, with Galois group generated
by the automorphism $\tau: y \mapsto y +1$ of order $p$. Let $P$ be a finite prime of
$k(x)$ and let $\nu_P$ be the corresponding valuation.  
Then $\nu_P(f(x)) \geq 0$, hence $P$ is unramified by \cite[Prop.\ III.7.8(b)]{sti}. 
For the infinite prime $\infty$ with corresponding valuation $\upsilon_\infty$, we have $$\nu_{{\infty}}(f(x)-(z^p-z)) \le 0$$ for all $z \in k[x]$ thus $P_{\infty}$ is totally ramified by \cite[Prop.\ III.7.8(c)]{sti}.

To prove part (2), we note that furthermore
$$\upsilon_{P_\infty}(y^p-y) = \upsilon_{P_\infty}(f(x))= \upsilon_{P_\infty}(x^s)  = -sp ,$$ which implies that $$\upsilon_{P_\infty}(y)=-s .$$ Now let $\widehat{\theta}$ be the completion of the valuation ring of $k(x)[y]/(y^p-y-f(x))$ at $P_{\infty}$, and let $\pi_\infty$ be a generator of the unique prime in $\widehat{\theta}$. 
Then write $y=\pi_\infty^{-s}u$, where $u$ is a unit in $\widehat{\theta} \simeq k[[\pi_\infty]]$. 
Since $k$ is algebraically closed, $\sqrt[s]{u} \in \widehat{\theta}$,  
and so $\sqrt[s]{y} \in \widehat{\theta}$.
After possibly changing $\pi_\infty$, 
we can assume without loss of generality that $\sqrt[s]{y} = \pi_\infty^{-1}$. 
Recalling that $\tau$ acts on $y$ by $\tau(y)=y+1$, we have
\begin{eqnarray}
\tau(\pi_\infty) & = & \tau\left(1/y \right)^{1/s}
= \left(\pi_\infty^s/(1+\pi_\infty^s) \right)^{1/s}
\nonumber\\ & = & \pi_\infty(1-\pi_\infty^s +
\pi_\infty^{2s}-+\ldots)^{\frac{1}{s}} \nonumber\\ & = &
\pi_\infty - (1/s)\pi_\infty^{s+1}+a_{2s+1}\pi_\infty^{2s+1}-+\ldots.
\nonumber
\end{eqnarray}
Thus $\upsilon_{P_\infty}(\tau(\pi_\infty) - \pi_\infty)=s+1$, which completes the proof of part (2).

To find the genus $g_Y$ of $Y$ for part (3), we make use of the Riemann-Hurwitz formula $$2g_Y-2=p(-2) + \sum_{i=0}^{\infty} \left(|I_i| -1\right), $$ where $I_i$ denotes the $i$th ramification group at $P_\infty$, \cite[Thms.\ 7.27 \& 11.72]{hkt}). From part (2), we then obtain that $g_Y
=(p-1)(s-1)/2$.
\end{proof}

Recall the following facts about the $p$th cyclotomic polynomial $\Phi_p(t):=t^{p-1}+\cdots +1$, 
which is the minimal polynomial over ${\mathbb Q}$ of a primitive $p$th root of unity $\zeta_p$. 
Now ${\mathbb Q}(\zeta_p)$ is a Galois extension of $\bbq$, unramified over $\ell$ since $\ell \not = p$,
and all primes over $\ell$ have the same residue field degree. 
The irreducible factors of $\Phi_p(t)$ modulo $\ell$ are in one-to-one correspondence with the primes of ${\mathbb Z}[\zeta_p]$ over $\ell$, 
and each of their degrees is equal to the residue field degree of the corresponding prime over $\ell$. 
The latter equals the order $a=\ord_p(\ell)$ of $\ell$ modulo $p$ \cite[Ch.\ 12.2, Exercise~\#20]{df}.

We shall soon explicitly construct a cover of $\bbp^1_k$ ramified only
over $\infty$ with Galois group $\Ga$.  
But before we do so, we start with a specific example.

\begin{exa}
{\rm Let $p$ be an odd prime. 
Consider the Artin-Schreier cover $\phi:Y_2 \to \mathbb{P}^1_k$ corresponding to the field extension 
$k(x) \hookrightarrow k(x)[y]/ (y^p - y - x^2)$. 
By Lemma~\ref{lem:genus}(3), the genus of $Y_2$ is $g_{Y}=(p-1)/2$.

Let $\Jac(Y)$ be the Jacobian of $Y$. 
The automorphism $\tau$ of $Y$ given by $\tau(y)=y+1$
defines an automorphism of $\Jac(Y)$ of order $p$. 

Now we describe the action of $\tau$ on the subgroup $\Jac(Y)[2]$ of 2-torsion points of $\Jac(Y)$ explicitly. 
Note that since $2g_Y=(p-1)$, then $\Jac(Y)[2]$ is isomorphic to $(\bbz/2\bbz)^{p-1}$
by \cite[pg.\ 64]{mumford}.
Thus we can represent $\tau$ as an element of $\GL_{p-1}(\bbz/2\bbz)$.

For $0 \leq i \leq p-1$, let $P_i$ denote the closed point of $Y$ at which the function $y-i$ vanishes.
For each $i$, the divisors $P_i$ and $D_i=P_i-P_\infty$ on $Y$ 
can be identified with elements of $\Jac(Y)$. 
Let $O$ be the identity element of $\Jac(Y)$, i.e., the linear equivalence class of principal divisors. 
Then the divisor $2D_i$ is equivalent to $O$ since ${\rm div}(y-i)=2 D_i$. 
Moreover since ${\rm div}(x)= D_0+D_1+ \cdots+D_{p-1}$ is equivalent to $0$, 
we have $D_i\in \Jac(Y)[2]$ with the only relation $D_{p-1}= -(D_0+D_1+\cdots+D_{p-2})$. 
In particular, $D_0, \ldots, D_{p-2}$ form a basis of $\Jac(Y)[2]$.
With respect to this basis, the action of $\tau$ can be represented by the $(p-1)\times (p-1)$-matrix
\begin{displaymath}
\left( \begin{array}{ccccc}
0 & 0  & \ldots & 0 & -1 \\
1 & 0  & \ldots & 0 & -1 \\
0 & 1  & \ldots & 0 & -1 \\
\vdots & \vdots & \ddots & 0 & -1 \\
0 & 0 & \ldots & 1 & -1 \\
\end{array} \right).
\end{displaymath}
The characteristic polynomial of $\tau$ is $\Phi_p(t)=1+t+ \ldots t^{p-1} \in (\bbz/2 \bbz)[t]$, 
which factors into irreducible polynomials each of degree equaling the order of $2$ modulo $p$. 
In particular, $\tau$ acts irreducibly on $\Jac(Y)[2]$ if and only if $2$ is a primitive root modulo $p$, 
i.e., if and only if $p$ is an Artin prime.

For example, if $p=3$, then $\tau$ acts irreducibly on $\Jac(Y)[2]$ 
with minimal polynomial $\Phi_3(t)=t^2+t+1$.
If $p=7$, then $2$ has order $3$ modulo $7$ and the factorization of $\Phi_7(t)$ 
into irreducible polynomials is
$\Phi_7(t) \equiv (x^3+x^2+1)(x^3+x+1)$ modulo $2$.
Thus the action of $\tau$ on $\Jac(Y)[2]$ can be represented by the $6 \times 6$-matrix
\begin{displaymath}
\left( \begin{array}{cc}
A_1 & 0 \\
0 & A_2 \\
\end{array} \right)
\end{displaymath}
where $A_1$ and $A_2$ are the 
irreducible $3$-dimensional companion matrices of $x^3+x^2+1$ and $x^3+x+1$ respectively.}
\end{exa}

For the rest of the paper, let $\phi_s: Y_s \to \mathbb{P}^1_k$ 
be the Artin-Schreier cover corresponding to the field extension $$k(x) \hookrightarrow k(x)[y]/ (y^p - y - x^s).$$
We show that $\phi_s$ can be dominated by a Galois cover of $\bbp^1_k$ with Galois group of the form $\Ga$
for $a$ equal to the order of $\ell$ modulo $p$.

\begin{prop} \label{Pexist}
Let $s$ and $\ell$  be primes distinct from $p$. Let $\phi_s:Y_s \to \bbp^1_k$ be the Artin-Schreier cover
with affine equation $y^p-y=x^s$. Let $a=\ord_p(\ell)$ be the order of $\ell$ modulo $p$. 
Then there exists an unramified Galois cover $\omega: Z_a \to Y_s$ with Galois group $\H^a$ 
such that $\psi_a=\phi_s \circ \omega: Z_a \rightarrow \bbp^1_k$ 
is a Galois cover of $\bbp^1_k$ ramified only over $\infty$ whose 
Galois group is a semi-direct product of the form $\H^a\rtimes \bbz/p\bbz$.
\end{prop}

\begin{proof}
By Lemma~\ref{lem:genus}(1), $\phi_s:Y_s \rightarrow \bbp^1_k$
is a Galois cover with Galois group $\bbz/ p \bbz$
ramified only at the point $P_\infty$ over $\infty$.
The genus $g_s$ of $Y_s$ is $(p-1)(s-1)/2$.
Consider two commuting automorphisms of $Y_s$ defined by
$$\tau: \begin{cases}
x \mapsto x,\\
y \mapsto y+1,
\end{cases} ~~ \sigma: \begin{cases}
x \mapsto \zeta_sx, \text{ where } \zeta_s \text{ is a primitive $s$th root of  unity,}\\
y \mapsto y.
\end{cases} $$

Let $\Jac(Y_s)$ be the Jacobian of $Y_s$. 
Then $\tau$ and $\sigma$ define commuting 
automorphisms of $\Jac(Y_s)$ of orders $p$ and $s$ respectively. 
Therefore, $\End(\Jac(Y_s))$ contains
a ring isomorphic to $\bbz[\zeta_p, \zeta_s] \cong \bbz[\zeta_{ps}]$, which is a $\bbz$-module of rank
$\phi(ps)=(p-1)(s-1)=2g_s$.
Then ${\mathbb Q}(\zeta_{ps})$ is contained in $\End(\Jac(Y_s)) \otimes {\mathbb Q}$. 
In other words, $\Jac(Y_s)$ has complex multiplication by ${\mathbb Q}(\zeta_{ps})$.

For a prime $\ell$ distinct from $p$, the automorphism $\tau$ induces an action on the subgroup $\Jac(Y_s)[\ell]$ of $\ell$-torsion points of $\Jac(Y_s)$. 
Recall that there is a bijection between $\ell$-torsion points $D$ of $\Jac(Y_s)$ and unramified $(\bbz/ \ell \bbz)$-Galois covers $\omega_D:Z_D \to Y_s$
\cite[Prop.\ 4.11]{milne}.  Also $D$ has order $\ell$ if and only if $Z_D$ is connected.
This induces a bijection between orbits of $\tau$ on the set of unramified $(\bbz/ \ell \bbz)$-Galois covers $\omega_D:Z_D \to Y_s$
and on the set of $\ell$-torsion points of $\Jac(Y_s)$. 
For a point $D$ of order $\ell$ of $\Jac(Y_s)$, consider the compositum $\omega: Z \to Y_s$ 
of all of the conjugates
$\omega_{\tau^j(D)}: Z_{\tau^i(D)} \to Y_s$ for $0 \leq j \leq p-1$: 
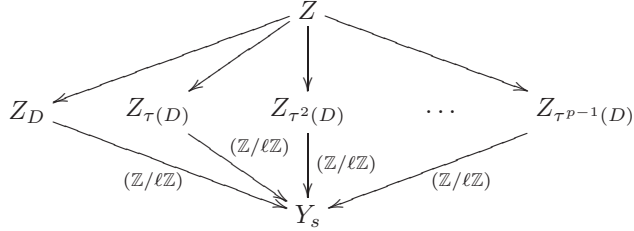
\begin{displaymath}
\xymatrix{
& & Z \ar[dll] \ar[dl] \ar[d] \ar[drr] & & \\
Z_D \ar[drr]_\H & Z_{\tau(D)} \ar[dr]^\H & Z_{\tau^2(D)} \ar[d]^\H & \ldots & Z_{\tau^{p-1}(D)} \ar[dll]^\H \\
& & Y_s & & }
\end{displaymath}

Then $Z$ is invariant under $\tau$ and so $\phi_s \circ \omega: Z \to \bbp^1_k$ is Galois.
Moreover, $\phi_s \circ \omega$ is the Galois closure of $\phi_s \circ \omega_D:Z_D \to \bbp^1_k$.

Suppose there is a non-trivial one-dimensional $\tau$-invariant subspace of $\Jac(Y_s)[\ell]$ with eigenvalue $1$; 
i.e. $\tau$ acts trivially on this subgroup of order $\ell$. 
This yields a cover $\psi_s \circ \omega_1:Z_1 \to Y_s \to \bbp^1_k$. 
Since the action of $\tau$ is trivial, $\psi_s \circ \omega_1$ is Galois, 
ramified only over $\infty$, with abelian Galois group 
$\bbz/\ell \bbz \times \bbz/p\bbz $. 
This contradicts Lemma \ref{Lgroup}. 

Since $\tau$ has order $p$, the minimal polynomial $m_\tau(t)$ of $\tau$ divides 
$t^p-1=(t-1)(t^{p-1}+\cdots +1)$ in $(\bbz/ \ell \bbz)[t]$. 
From the preceding paragraph, there is no non-trivial one-dimensional $\tau$-invariant subspace of 
$\Jac(Y_s)[\ell]$ with eigenvalue $1$.
This implies that $m_\tau(t)$ divides the $p$th cyclotomic polynomial 
$\Phi_p(t)=t^{p-1}+\cdots +1$ in $(\bbz/ \ell \bbz)[t]$. 
The irreducible factors of $\Phi_p(t)$ in $(\bbz/ \ell \bbz)[t]$ all have degree $a$.
Thus the degree of $m_\tau(t)$ equals $a$.
 
Since $2g_s=(p-1)(s-1)$, we have $\Jac(Y_s)[\ell] \cong (\bbz/\ell\bbz)^{(p-1)(s-1)}$, 
so we can represent $\tau$ as an element of 
$\GL_{(p-1)(s-1)}(\bbz/\ell\bbz)$.  
We can choose a basis of $\Jac(Y_s)[\ell]$ such that $\tau$ is represented as an 
element of $\GL_{(p-1)(s-1)}(\bbz/\ell\bbz)$ in block form.
The first irreducible subrepresentation of $\tau$ has dimension $a$. 
Moreover, since ${\mathbb Q}(\zeta_{ps})$ is a Galois extension of ${\mathbb Q}$, the block form of $\tau$ consists entirely of irreducible blocks of the same size. 
In particular, the number of irreducible blocks is $(p-1)(s-1)/a$.  In other words, $\tau$ can be represented by an element of $\GL_{(s-1)(p-1)}(\bbz/\ell\bbz)$
of the form 

\begin{displaymath}
\left( \begin{array}{ccc}
A_1 & & 0 \\
& A_2 &  \\
& \ddots &  \\
0 & & A_{(p-1)(s-1)/a} 
\end{array} \right), 
\end{displaymath}
where $A_i$ is an $a\times a$ matrix representing an $a$-dimensional irreducible subrepresentation of $\tau$ on $\Jac(Y_s)[\ell]$.

Using the bijection between orbits of $\Jac(Y_s)[\ell]$ and orbits of $(\bbz/\ell\bbz)$-covers of $Y_s$ 
under $\tau$ and the above observation for the action of $\tau$ on $\Jac(Y_s)[\ell]$, there exists an
unramified $(\bbz/\ell \bbz)^a$-Galois cover $\omega: Z_a \to Y_s$
such that $\psi_a=\phi_s \circ \omega: Z_a \rightarrow \bbp^1_k$ 
is a Galois cover of $\bbp^1_k$ with Galois group of the form $\H^a\rtimes \bbz/p\bbz$.
Also $\psi_a$ is ramified only over infinity 
since $\phi_s$ is ramified only over $\infty$ and since $\omega$ is unramified.
\end{proof}

\section{Minimal genus of $\Gb$-Galois covers of $\bba^1_k$}\label{unique}

In this section, we find the minimal genus of a curve $Z$ that admits a covering map 
$\psi: Z \to \bbp^1_k$ ramified only over $\infty$, with Galois group of the form $\Gb$.
The minimal genus depends only on $\ell$ and $p$.
We consider the cases $p$ odd and $p=2$ separately.
We also prove that the number of curves $Z$ of this minimal genus which admit such a covering map
is at most $(p-1)/a$ when $p$ is odd and at most $\ell+1$ when $p=2$.
The following lemma will be useful.

\begin{lem} \label{Lunram}
Let $G$ be a semi-direct product of the form $\Gb$ where $\ell$ is a prime distinct from $p$.
If $\psi:Z \to \bbp^1_k$ is a Galois cover ramified only over $\infty$ with Galois group $G$, 
then the subcover $\omega:Z \to Y$ with Galois group $\H^b$ is unramified.
\end{lem}

\begin{proof}
The quotient of $G$ by the normal subgroup $N=({\mathbb Z}/\ell {\mathbb Z})^b$ 
is ${\mathbb Z}/p {\mathbb Z}$.  
Thus the cover $\psi$ is a composition $\psi=\phi \circ \omega$ where 
$\phi:Y \to \bbp^1_k$ has Galois group ${\mathbb Z}/p {\mathbb Z}$ and 
is totally ramified at the unique point $P_\infty$ over $\infty$
and where $\omega:Z \to Y$ has Galois group $N$ and is branched only over $P_\infty$.
Then $\omega$ is a prime-to-$p$ abelian cover of $Y$.  
Let $g$ be the genus of $Y$.  
Then by \cite[XIII, Cor.\ 2.12]{GR1}, the prime-to-$p$ fundamental group of $Y-\{P_\infty\}$ 
is isomorphic to the prime-to-$p$ quotient $\Gamma$ of the 
free group on generators $\{a_1,b_1, \ldots, a_g,b_g, c\}$ subject to the relation 
$\prod_{i=1}^{g}[a_i,b_i] = c^{-1}$.  
The cover $\omega$ corresponds to a surjection of $\Gamma$ onto $N$ 
where $c$ maps to the canonical generator of inertia $\gamma$ of a point of $Z$ over $P_\infty$.  
Thus $N$ is generated by elements  $\{\alpha_1,\beta_1, \ldots, \alpha_g,\beta_g, \gamma \}$ subject to the relation $\prod_{i=1}^{g}[\alpha_i,\beta_i] = \gamma^{-1}$.
Then $\gamma=1$ since $N$ is abelian and so $\omega$ is unramified.
\end{proof}

\begin{thm}\label{exist-minimal}
Let $p$ be an odd prime.  Let $\ell$ be a prime distinct from $p$ and 
let $a$ be the order of $\ell$ modulo $p$. Then:
\begin{enumerate}
\item There exists a Galois cover $\psi_a:Z_a \to \bbp^1_k$ ramified only over $\infty$ 
whose Galois group is a semi-direct product of the form $\Ga$ such that $g_{Z_a}=1+ \ell^a(p-3)/2$.
\item The integer $g_{Z_a}$ is the minimal genus of a curve $Z$ which admits a covering map $\psi:Z \to \bbp^1_k$
ramified only over $\infty$ with Galois group of the form $\Gb$ for any positive integer $b$.
\item There are at most $(p-1)/a$ isomorphism classes of curves $Z$ 
which admit a Galois covering map as in part (1) with minimal genus $g_{Z_a}$.
\end{enumerate}
\end{thm}

\begin{proof}
By the construction in Proposition \ref{Pexist}, 
there exists a Galois cover $\psi_a:Z_a \to \bbp^1_k$ ramified only over $\infty$ 
whose Galois group is a semi-direct product of the form $\Ga$. 
We compute the genus of the curve $Z_a$.
Recall that $\psi_a$ is a composition $\psi=\phi_2 \circ \omega$ where 
$\omega:Z \to Y_2$ is an unramified $\H^a$-Galois cover and 
$\phi_2: Y_2 \to \bbp^1_k$ has Artin-Schreier equation $y^p-y=x^2$.
Then $Y_2$ has genus $g_{Y_2}=(p-1)/2$ by Lemma \ref{lem:genus}(3).
By the Riemann-Hurwitz formula, $2g_{Z_a}-2=\ell^a(2g_{Y_2}-2)=\ell^a(p-3)$, 
i.e., $g_{Z_a}=1+\ell^a(p-3)/2$.  This completes part (1).

For part (2), 
suppose $\psi:Z \to \bbp^1_k$ is a Galois cover ramified only over $\infty$ with Galois group of the form $\Gb$. 
If $g$ is the genus of $Z$, we will show that $g \geq g_{Z_a}$.
As described in the proof of Lemma \ref{Lunram}, 
the cover $\psi$ is a composition $\psi=\phi \circ \omega$ where  
$\phi: Y \to \bbp^1_k$ has Galois group $\bbz /p \bbz$ and is ramified only over $\infty$ 
and where $\omega$ is unramified with group $\H^b$.
By the Riemann-Hurwitz formula, $2g-2=\ell^b(2g_Y-2)$.

By Artin-Schreier theory, $\phi$ is given by an equation $y^p-y=f(x)$ where $f \in k[x]$ has degree $s$ for some integer $s$ relatively prime to $p$.
Since the genus $g_Y$ of $Y$ is $(p-1)(s-1)/2$ by Lemma~\ref{lem:genus}~(3), 
we should make $s$ as small as possible. 
The value $s=1$ is impossible since then $Y$ is a projective line and 
there do not exist Galois covers of the projective line 
ramified only over one point with Galois group $\bbz/\ell \bbz$. 
Thus $s=2$ yields the smallest possible value for $g_Y$, namely $(p-1)/2$.
Recall that $b \geq a$ by Lemma \ref{Lgroup}.
Thus $g \geq 1 + \ell^a(p-3)/2 = g_{Z_a}$.

For part (3), suppose $\psi:Z \to \bbp^1_k$ is a Galois cover ramified only over $\infty$ with Galois group of the form $\Ga$ and the genus of $Z$ satisfies $g_Z=1+\ell^a(p-3)/2$.  
As in part (2), $\psi$ factors as $\phi \circ \omega$ where $\omega:Z \to Y$ is an unramified $\H^a$-Galois cover, where
$\phi:Y \to \bbp^1_k$ is an Artin-Schreier cover ramified only over $\infty$, and where $Y$ has genus $(p-1)/2$.
By Lemma~\ref{lem:genus}(3), $Y$ has an affine equation $y^p-y=a_2x^2+a_1x+a_0$ for some $a_0,a_1 \in k$ and $a_2 \in k^*$.
Since $p$ is odd and $k$ is algebraically closed, it is possible to complete the square
and write $a_2x^2+a_1x+a_0=x_1^2+\epsilon$ where $x_1=\sqrt{a_2}x+a_1/2\sqrt{a_2}$.
After modifying by an automorphism of the projective line, specifically by the affine 
linear transformation $x \mapsto x_1$,
the equation for $Y$ can be rewritten as $y^p-y=x_1^2+\epsilon$.
Since $k$ is algebraically closed, there exists $\delta \in k$ such that $\delta^p-\delta=\epsilon$.  Let $y_1=y-\delta$.
After the change of variables $y \mapsto y_1$, the curve $Y$ is isomorphic to the curve $Y_2$ with affine equation $y_1^p-y_1=x_1^2$.
Thus there is a unique possibility for the isomorphism class of the curve $Y$.

From the proof of Proposition~\ref{Pexist}, there is a bijection between 
$\tau$-invariant connected unramified $\H^a$-Galois covers of $Y_2$ and 
orbits of $\tau$ on points $D$ of order $\ell$ on $\Jac(Y_2)$.
The action of $\tau$ on $\Jac(Y_2)[\ell]$ decomposes into $(p-1)/a$ irreducible subrepresentations.
Each of these is distinct, because the irreducible 
factors of $\Phi_p(t)\in (\bbz/\ell \bbz)[t]$ are distinct.
Thus there are $(p-1)/a$ choices for a $\tau$-invariant unramified $\H^a$-Galois cover of $Y_2$. 
Thus there are at most $(p-1)/a$ isomorphism classes of curves $Z$ 
which admit a Galois covering map as in part (1) with minimal genus $g_{Z_a}$.
\end{proof}

We note that the set of curves which are unramified $\H^a$-Galois covers of $Y_2$ 
may contain fewer than $(p-1)/a$ isomorphism classes of curves.

\begin{thm}\label{p=2exist-minimal}
Let $p=2$ and let $\ell$ be an odd prime.  Then:
\begin{enumerate}
\item There exists a Galois cover $\psi:Z \to \bbp^1_k$ ramified only over $\infty$ 
with Galois group of the form $\bbz/\ell \bbz \rtimes \bbz/2 \bbz$.
\item The minimal genus of a curve $Z$ which admits a covering map as in part (1) is $g_Z=1$.
\item There are at most $\ell + 1$ isomorphism classes of curves $Z$ 
which admit a Galois covering map as in part (1) with minimal genus $g_Z=1$.
\end{enumerate}
\end{thm}

\begin{proof}
Note that the order of $\ell$ modulo $2$ is $a=1$.
For part (1), Lemma~\ref{LexistGb} shows that there exists a semi-direct product of the form 
$\bbz/\ell \bbz \rtimes \bbz/2 \bbz$ which is quasi-$2$. 
The result is then immediate from Theorem~\ref{Tabconj}.

Suppose $\psi:Z \to \bbp^1_k$ is a Galois cover ramified only over $\infty$ with Galois group 
as in part (1).
As before, $\psi$ factors as a composition $\phi \circ \omega$.
where $\omega:Z \to Y$ has Galois group $\bbz/\ell \bbz$ 
and $\phi:Y \to \bbp^1_k$ is an Artin-Schreier extension 
with affine equation $y^2-y=f(x)$ for some $f(x) \in k[x]$ of odd degree $s$.
By Lemma \ref{Lunram}, $\omega$ is unramified.
The minimal genus for $Z$ will thus occur when $s$ is as small as possible.
As before, $s=1$ is impossible, and so $s=3$ is the smallest choice.
In this case, by Lemma~\ref{lem:genus}(3), $g_Y=1$, i.e., $Y$ is an elliptic curve. 
By the Riemann-Hurwitz formula, the minimal genus for $Z$ is $g_Z=1 + \ell(g_Y-1)=1$, which completes part (2).

For part (3), since $k$ is algebraically closed, we can 
complete the cube of $f(x)$ and make the corresponding change of variables, which is 
a scaling and translation of $x$.
So we can assume that $Y$ has affine equation $y^2-y=x^3+a_1x+a_0$ 
for some $a_0, a_1 \in k$. 
Then it follows from \cite[Appendix A,~Prop.\ 1.1c]{sil} that the $j$-invariant of $Y$ is $j(Y)=0$ and that the discriminant is $\Delta(Y)=(-1)^4=1$.  
Since $k$ is algebraically closed, by \cite[Appendix A,~Prop.\ 1.2b]{sil}, all elliptic curves $Y$ with $j(Y)=0$ are isomorphic over $k$.
Thus there is a unique choice for $Y$ up to isomorphism. 
Without loss of generality, we may assume that $Y=Y_3$ has affine equation $y^2-y=x^3$.

From the proof of Proposition~\ref{Pexist}, the action of $\tau$ on $\Jac(Y_3)[\ell]$ decomposes into the direct sum of two $1$-dimensional subrepresentations.  In other words, the action of $\tau$ is diagonal with both eigenvalues equal to $-1$.
The number of non-trivial $\tau$-invariant subgroups of $\Jac(Y_3)[\ell]$ is the number of subgroups of 
order $\ell$ in $(\bbz/\ell \bbz)^2$, which is $\ell + 1$.
As in Theorem \ref{exist-minimal}, this implies that there are at most $\ell +1$ isomorphism classes of curves $Z$ 
which admit a Galois covering map as in part (1) with minimal genus $g_Z=1$.
\end{proof}

We note that the set of curves which are unramified $\bbz/\ell \bbz$-Galois covers of $Y_3$ 
may contain fewer than $\ell +1$ isomorphism classes of curves.

\end{document}